\documentclass[final]{svjour3}

\usepackage{amsmath}
\usepackage{amssymb}
\usepackage{amsfonts}
\usepackage{graphicx}
\usepackage{psfrag}
\usepackage{url}
\usepackage{verbatim}
\usepackage{cite}
\usepackage{subfig}
\usepackage{tikz}
\usetikzlibrary{positioning,fit,arrows,calc}
\usepackage[fleqn,tbtags]{mathtools}

\def\MATLAB{{\sc Matlab}}

\def\E{{\rm e}}
\def\D{\,\hbox{d}}

\def\sech{{\rm sech}}

\title{Fast inverse transform sampling in one and two dimensions}

\author{Sheehan Olver \and Alex Townsend}

\institute{School of Mathematics and Statistics, The University of Sydney, NSW 2006, Australia. (sheehan.olver@sydney.edu.au, http://www.maths.usyd.edu.au/u/olver/). Supported by the ARC grant DE130100333 \and Mathematical Institute, 24-29 St Giles', Oxford,  OX1 3LB, UK. (townsend@maths.ox.ac.uk, http://people.maths.ox.ac.uk/townsend/). Supported by the EPSRC grant EP/P505666/1 and the ERC grant FP7/2007-2013 to Nick Trefethen} 

\begin{document}
\maketitle

\begin{abstract}
	We develop a computationally efficient and robust algorithm for generating
pseudo-random samples from a broad class of smooth probability distributions in one and two dimensions.   The  algorithm is based on inverse transform sampling with a polynomial approximation
scheme using Chebyshev polynomials, Chebyshev grids, and low rank function approximation.
Numerical experiments demonstrate that our algorithm outperforms existing approaches.
\end{abstract}

\section{Introduction}

Generating pseudo-random samples from a prescribed probability distribution is an important task in applied probability, statistics,  computing, physics, and elsewhere. A classical approach
is \emph{inverse transform sampling}, in which pseudo-random samples $U_1,\ldots,U_N$ are generated from a uniform distribution 
$U$ on $[0,1]$ and then transformed by $F_X^{-1}(U_1),\ldots,F_X^{-1}(U_N)$, where $F_X$ is the cumulative distribution function  (CDF) of a prescribed probability distribution. 

Conventional wisdom is that the practical application of inverse transform sampling is limited.
 To quote \cite[pp. 103]{GentleSampling}:
\begin{quote}
[Inverse transform sampling] requires a complete approximation to [the CDF] regardless of the desired sample size, it does not generalize to multiple dimensions and it is less efficient than other approaches.
\end{quote}
Similar comments are present in other computational statistics texts, for example, \cite[pp. 124]{StatisticalAtmosphere}, \cite[pp. 155]{GivensComputationalStats}, and \cite[pp. 122]{TanizakiComputational}.  
Instead, other methods are advocated for low dimensions, such as slice sampling \cite{Neal_03_01} or rejection sampling \cite{TanizakiComputational,StatisticalAtmosphere,GentleSampling,GivensComputationalStats}.  In higher dimensions,  algorithms concentrate on avoiding the curse of dimensionality, such as Metropolis--Hastings sampling \cite{MetropolisSampling}.
Moreover, there are many other methods that are designed to exploit certain special properties: e.g.,  ziggurat algorithm \cite{Ziggurat} (requires monotonicity) and adaptive rejection sampling \cite{Gilks_92_01,Gilks_92_02} (requires log-convexity).  

In contradiction to conventional wisdom,  we show that inverse transform sampling is computationally efficient and robust  when combined with an adaptive Chebyshev polynomial approximation scheme.  Approximation by Chebyshev polynomials converges superalgebraically fast for smooth functions, and hence the CDF can be numerically computed to very high accuracy, i.e., close to machine precision.   Rather than forming the inverse CDF --- which is often nearly singular --- we invert the CDF pointwise using a simple bijection method.   This approach is robust --- it is guaranteed to converge to a high accuracy --- and fast, due to the simple representation of the CDF.

Furthermore, we extend our approach to probability distributions of two 
variables by exploiting recent developments in low rank function approximation
to achieve a highly efficient inverse transform sampling scheme in 2D.   
Numerical experiments confirm that this approach outperforms existing sampling 
routines for many black box distributions, see Section~\ref{sec:numericalexamples}.

\begin{remark}
We use {\it black box distribution} to mean a  continuous probability distribution that can only be evaluated pointwise\footnote{Strictly, we require the distribution 
to be continuous with bounded variation \cite{Trefethen_13_01}, and in practice this is almost always satisfied by continuous functions.  We do not require the distribution to be prescaled to integrate to one, nor that the distribution decays smoothly to zero.}.  The only other piece of information we assume is 
an interval $[a,b]$, or a rectangular domain $[a,b]\times [c,d]$, containing the numerical support of the distribution.  
\end{remark}



A \MATLAB\ implementation is publicly available from \cite{InverseTransformCode}.  As an example, the following syntax generates $2,\!000$ samples from a 2D distribution: 
	\begin{verbatim}
		[X,Y] = sample(@(x,y) exp(-x.^2-2.*y.^2).*(x-y).^2,[-3,3],[-3,3],2000);
	\end{verbatim}
The algorithm scales the input function to integrate to one, and generates pseudo-random samples from the resulting probability distribution.


%

\section{Inverse transform sampling with polynomial approximation}\label{sec:inversesampling}

Let $f(x)$ be a prescribed non-negative function supported on the interval $[a,b]$.  Otherwise, if the support of $f$ is the whole real line and $f$ is rapidly decaying, an interval $[a,b]$ can be selected outside of which $f$ is negligible.  The algorithm will draw samples that match the probability distribution
	$$f_X(x) = {f(x) \over \int_a^b f(s) \D s}$$
up to a desired tolerance $\epsilon$ (typically close to machine precision such as $2.2\times10^{-16}$).  

	First, we replace $f$ by a polynomial approximant $\tilde{f}$ on $[a,b]$, which we numerically compute such that $\|f - \tilde{f}\|_\infty \leq \epsilon\|f\|_\infty$.  That is, we 
construct a polynomial approximant $\tilde{f}$ of $f$ in the form
\begin{equation}
\tilde{f}(x) = \sum_{k=0}^n \alpha_k T_k\left(\frac{2(x-a)}{b-a}-1\right), \quad \alpha_k\in\mathbb{R}, \qquad x\in[a,b], 
\label{eq:ChebExpansion}
\end{equation}
where $T_k(x) = \cos(k\cos^{-1}x)$ is the degree $k$ Chebyshev polynomial. While the polynomial approximants \eqref{eq:ChebExpansion} could be represented in other polynomial bases --- such as the monomial basis ---  the Chebyshev basis is a particularly good choice for numerical computation as the expansion coefficients $\alpha_0,\ldots,\alpha_n$ can be stably calculated in $\mathcal{O}(n\log n)$ operations \cite{Mason_03_01,Trefethen_13_01}.  This is accomplished by applying the fast cosine transform to function evaluations of $f$ taken from a Chebyshev grid in $[a,b]$ 
of size $n=2^k+1$, i.e., at the set of points
\[
\left\{(a+b)/2 + (b-a)/2\cos(j\pi/n)\right\}_{0\leq j\leq n}.
\]
We can adaptively select $n$ by 
sampling $f$ at Chebyshev grids of size $2^3 + 1 \,(= 9)$, $2^4+1 \,(= 17)$, $2^5 + 1\,(= 33)$, and so on, 
until the tail of the coefficients $\left\{\alpha_j\right\}_{0\leq j\leq n}$ has decayed to below machine precision relative to the absolute maximum of $f$.     In addition to the fast transform, Chebyshev expansions are numerically stable \cite{Trefethen_13_01},  in contrast to the 
well-known numerical instability associated to interpolation at equispaced points \cite{Platte_11_01}.



Instead of symbolically manipulating $f$, which is not possible if $f$ is a black box distribution, we 
numerically manipulate $\tilde{f}$ to apply the inverse transform sampling approach.   This rests on fast and well-established algorithms 
for evaluating and integrating polynomial interpolants of the form \eqref{eq:ChebExpansion}.  (A convenient implementation of many of these algorithms 
is in the {\sc Chebfun} software system \cite{Chebfun}.)     Figure~\ref{fig:InverseSamplingChebyshev} summarizes the approach we are about to describe based 
on the inverse transform sampling and Chebyshev approximation.

\begin{figure}
\centering
	\begin{tabular}{|l|}
	\hline
	\rule{0pt}{3ex}
	\textbf{Algorithm: Inverse transform sampling with Chebyshev approximation}\\[10pt]
	\textbf{Input:} A non-negative function $f(x)$, a finite interval $[a,b]$, and an integer $N$.\\[3pt]
	\textbf{Output:} Random samples $X_1,\ldots,X_N$ drawn from $f_X = f/\int_a^b f(s) \D s$.\\[5pt]
	Construct an approximant $\tilde{f}$ of $f$ on $[a,b]$. \\[3pt]
        Compute $\tilde{F}(x) = \int_a^x \tilde{f}(s) \D s$ using \eqref{eq:cumsumrelation}.\\[3pt]
        Scale to a CDF, $\tilde{F}_X(x) = \tilde{F}(x)/\tilde{F}(b)$, if necessary. \\[3pt]
        Generate random samples $U_1,\ldots,U_N$ from the uniform distribution on $[0,1]$.\\[3pt]
	\textbf{for $k=1,\ldots,N$ } \\[5pt]
	\hspace{.75cm} Solve $F_X(X_k) = U_k$ for $X_k$.\\[4pt]
        \textbf{end}\\[3pt]
	\hline
	\end{tabular}
	\caption{Pseudocode for the inverse transform sampling algorithm with polynomial approximation using Chebyshev polynomials and Chebyshev grids. 
}
	\label{fig:InverseSamplingChebyshev}
\end{figure}

With the polynomial approximation $\tilde f$ in hand, we  construct the corresponding cumulative function, denoted by $\tilde{F}$, by calculating the indefinite integral 
\begin{align*}
\tilde{F}(x) &= \int_a^x \tilde{f}(s) \D s = \sum_{k=0}^n \alpha_k \int_{a}^x T_k\left(\frac{2(s-a)}{b-a}-1\right) \D s \\
&= \sum_{k=0}^n \alpha_k \frac{b-a}{2} \int_{-1}^{\frac{2(x-a)}{b-a}-1} T_k(t) \D t,
\end{align*}
where the last equality comes from a change of variables with $t = \frac{2(s-a)}{b-a}-1$. Fortunately, we have the following relation satisfied by 
Chebyshev polynomials \cite[pp.~32--33]{Mason_03_01}:
\begin{equation}
\int^{s} T_k(t) \D t = \begin{cases} 
                          \frac{1}{2}\left(\frac{T_{k+1}(s)}{k+1} - \frac{T_{|k-1|}(s)}{k-1}\right), & k\not = 1, \\[3pt]
	                  \frac{1}{4}T_{2}(s), & k=1, \\[3pt]
                          \end{cases}
\label{eq:cumsumrelation}
\end{equation}
where $s\in[-1,1]$. Therefore, using \eqref{eq:cumsumrelation}, we can
compute coefficients $\beta_0,\ldots,\beta_n$ in $\mathcal{O}(n)$ operations (see, for example, \cite[pp.~32--33]{Mason_03_01}) such that 
\begin{equation}
\tilde{F}(x) = \sum_{k=0}^n \beta_k T_k\left(2\frac{x-a}{b-a}-1\right), \quad \beta_k\in\mathbb{R}, \qquad x\in[a,b]. 
\label{eq:Chebcumsum}
\end{equation} 
At this stage we
rescale  $\tilde{F}$ so that it is a CDF: $\tilde{F}_X(x) = \tilde F(x) / \tilde F(b)$.

Once we have the approximant $\tilde{F}_X$, given in \eqref{eq:Chebcumsum}, we generate $N$ pseudo-random samples $U_1,\ldots, U_N$ from the uniform distribution on $[0,1]$ using 
a standard inbuilt pseudo-random number generator, and then solve the following 
rootfinding problems for $X_1,\ldots,X_N$:
\begin{equation}
\tilde{F}_X(X_k) = U_k,\qquad k=1,\ldots,N.
\label{eq:rootfinding}
\end{equation}
The rootfinding problems in \eqref{eq:rootfinding} can be solved by various standard methods such as 
the eigenvalues of a certain matrix \cite{Trefethen_13_01}, Newton's method, the bisection method, or Brent's method \cite[Chapter 4]{Brent_73_01}.

To achieve robustness --- guaranteed convergence, computation cost, and accuracy --- we  use the bisection method.  
The bisection method approximates the (unique) solution $X_k$ to $F_X(X_k)= U_k$ by a sequence of intervals $[a_j,b_j]$ that contain $X_k$, satisfying $F_X(a_j) < U_k < F_X(b_j)$.   The initial interval is the whole domain $[a,b]$.  Each stage finds the midpoint $m_j$ of the interval $[a_j,b_j]$, choosing the new interval based on whether $F_X(m_j)$ is greater or less than $U_k$.  Convergence occurs when ${b_j - a_j} < tol (b - a)$, where $b_j-a_j = 2^{-j}(b-a)$.  For example,  it takes precisely 47 iterations to converge when $tol = 10^{-14}$ and a further five iterations to converge to machine precision. Since the CDF is represented as a Chebyshev series \eqref{eq:Chebcumsum} it can be efficiently evaluated using Clenshaw's method, which is an extension of Horner's scheme to Chebyshev series \cite{Clenshaw_55_01}.

The inverse transform sampling with Chebyshev approximation is very efficient, as demonstrated  in the numerical experiments of Section \ref{sec:numericalexamples} using the \MATLAB\ implementation \cite{InverseTransformCode}. 


\section{Inverse transform sampling in two variables}\label{sec:twodimensions}

We now extend the inverse transform sampling approach with Chebyshev approximation to probability distributions of two variables.
Figure~\ref{fig:InverseSamplingChebyshev2} summarizes our algorithm, which generates $N$ pseudo-random samples $(X_1,Y_1),\ldots,(X_N,Y_N)$ drawn 
from the probability distribution $f_{X,Y}(x,y) = f(x,y)/\int_a^b \int_c^d f(s,t) \D s \D t$.   The essential idea is to replace $f$ by a low rank approximant $\tilde f$ that can be manipulated in a computationally efficient manner.  
Again, for convenience, if $f$ is a non-negative function that does not integrate to one then we automatically scale the function to a 
probability distribution. 


\begin{figure}
\centering
	\begin{tabular}{|l|}
	\hline
	\rule{0pt}{3ex}
	\textbf{Algorithm: Inverse transform sampling for $f_{X,Y}$}\\[10pt]
	\textbf{Input:} A non-negative function $f(x,y)$,  a finite domain $[a,b]\times [c,d]$, and an integer $N$.\\[3pt]
	\textbf{Output:} Random samples $(X_1,Y_1),\ldots,(X_N,Y_N)$ from $f_{X,Y} = f/\int_c^d\int_a^b f(s,t)\D s \D t$.\\[5pt]
	Construct a low rank approximant $\tilde{f}$ of $f$ on $[a,b]\times [c,d]$. \\[3pt]
        Compute $\tilde{f}_1(x) = \int_c^d \tilde{f}(x,y) \D y$. \\[3pt]
        Generate random samples $X_1,\ldots,X_N$ from $f_X(x) = \tilde{f}_1(x)/\tilde f_1(b)$.\\[3pt]
        
	\textbf{for $i=1,\ldots,N$ } \\[5pt]
	\hspace{.75cm} Generate a random sample $Y_k$ from $\tilde{f}_{Y|X=X_i}(y) = \tilde{f}(X_k,y)/\tilde f_1(X_k)$.\\[4pt]
        \textbf{end}\\[3pt]
	\hline
	\end{tabular}
	\caption{Pseudocode for the inverse transform sampling algorithm for a probability distribution of two variables denoted by $f_{X,Y}$. 
}
	\label{fig:InverseSamplingChebyshev2}
\end{figure}

We use a {\it low rank approximation} of $f(x,y)$.  
A rank-$1$ function is a function of two variables that can be written as the product of two univariate functions, i.e., $r_1(x) c_1(y)$. Moreover, a rank-$k$ function 
is a sum of $k$ rank-1 functions. We efficiently approximate a
probability distribution $f$ by a rank-$k$ function by using Gaussian elimination \cite{Townsend_13_01}:
\begin{equation}
f(x,y) \approx \tilde{f}(x,y) = \sum_{j=1}^k \sigma_j r_j(x)c_j(y),
\label{eq:lowrank}
\end{equation}
where $k$ is adaptively chosen and $c_j$ and $r_j$ are polynomials of one variable represented in Chebyshev expansions of the form \eqref{eq:ChebExpansion}.   For simplicity we assume that
$c_j$ and $r_j$ are polynomials of degree $m$ and $n$, respectively, and then a function is considered to be of low rank if $k\ll \min(m,n)$. 

While most functions are mathematically of infinite rank  --- such as $\cos(xy)$ ---  they
can typically  be  approximated by low rank functions to high accuracy, especially if the function is smooth \cite{Townsend_13_02}. It is well-known that the singular value decomposition can be used to compute optimal low rank approximations of matrices, and it can easily be extended for the approximation of functions. 
However, the Gaussian elimination approach from \cite{Townsend_13_01}  is significantly faster and  constructs near-optimal low rank function approximations \cite{Bebendorf_08_01,Townsend_13_01}.  This algorithm is conveniently implemented in the recent software package {\sc Chebfun2} \cite{Townsend_13_02}, which we utilize to construct  $\tilde{f}$. Furthermore, we emphasize that this low rank approximation process only requires pointwise evaluations of $f$. 

With the low rank approximation $\tilde f$, given in \eqref{eq:lowrank}, we can efficiently perform the steps of the 2D inverse transform sampling algorithm.
First,  we approximate the marginal distribution of $X$ by integrating $\tilde{f}$ over the second variable:
\[
\tilde{f}_1(x)  = \int_c^d \tilde{f}(x,y) \D y = \sum_{j=1}^k\sigma_j r_j(x) \int_c^d c_j(y)\D y.
\]
Therefore, $\tilde{f}_1$ can be computed by applying a 1D quadrature rule, such as the Clenshaw--Curtis quadrature rule \cite{Clenshaw_60_01}, to $k$ polynomials of one variable (rather than anything inherently 2D). The resulting sum 
is  a Chebyshev expansion (in $x$), and we use the 
algorithm described in Section~\ref{sec:inversesampling} to generate $N$ pseudo-random samples $X_1,\ldots, X_N$ from $\tilde{f}_X(x) = \tilde f_1(x)/\tilde f_1(b)$.  

Afterwards, we generate the corresponding $Y_1,\ldots,Y_N$ by using a numerical approximation to the conditional distribution denoted by $\tilde{f}_{Y|X}$. To construct $\tilde{f}_{Y|X}$ for
each $X_1,\ldots, X_N$ we require evaluation, i.e.,  
\[
\tilde{f}_{Y|X=X_i}(y) =  \frac{\tilde{f}(X_i,y)}{\tilde{f}_1(X_i)}, \qquad i=1,2,\ldots,N,
\]
where $\tilde f_1(X_i) = \int_{c}^d \tilde{f}(X_k,t) \D t$ is precisely the normalization constant.  
Fortunately, evaluation of $\tilde{f}(X_i,y)$ is relatively cheap because of the low rank representation of $\tilde{f}$. That is,
\[
\tilde{f}(X_i,y) = \sum_{j=1}^k \sigma_j r_j(X_i)c_j(y),
\]
and thus we only require the evaluation of $k$ polynomials of one variable (again, nothing inherently 2D), which, as before, is accomplished using Clenshaw's algorithm \cite{Clenshaw_55_01}. The final task of generating the sample $Y_i$ has been reduced to the 
1D problem of drawing a sample from $\tilde f_{Y|X=X_i}$, which was solved in Section~\ref{sec:inversesampling}.

In total we have generated $N$ pseudo-random samples $(X_1,Y_1),\ldots, (X_N,Y_N)$ from the probability distribution $f_{X,Y}$ without an explicit expression for the CDF or its inverse. Low rank approximation 
was important for computational efficiency, and the total computational cost for $N$ samples, including the low rank approximation of $f$, is
	\begin{align*}
	\overbrace{\mathcal{O}(k^3 + k^2(m+n))}^{\hbox{low rank approximation}}+ &  \overbrace{\mathcal{O}(k(m+n))}^{\hbox{marginal distribution}} + \overbrace{\mathcal{O}(Nn)}^{\hbox{sample $X_i$}}  \cr
	& + \overbrace{\mathcal{O}(N (m + kn))}^{\hbox{construct $\tilde f_{Y|X = X_i}(y)$}} + \overbrace{\mathcal{O}(N  m)}^{\hbox{sample $Y_i$}} \cr
	= \mathcal{O}(k^3 + k^2 (m+n) + N ( n+ k m) ).
	\end{align*}
Therefore, provided the probability distribution can be well-approximated by a low rank function, the inverse transform sampling approach in 2D has the same order of complexity --- i.e., linear in $N$, $m$ and $n$ --- as the 1D approach described in Section~\ref{sec:inversesampling}.  Many probability distribution can be approximated by functions of low rank, as we see below.  

\section{Numerical experiments}\label{sec:numericalexamples}

\begin{table}
\begin{center}
\begin{tabular}{c | c || c | c | c  c |  c |  c} 
		PDF 										& Rank 	& ITS  		&  SS  		& RS \\
		\hline 
		$\E^{-{x^2 \over 2}} (1+ \sin^2 3 x)(1+ \cos^2(5 x))$ 	&---		&  0.29 		&   3.01	  	& 0.55 \cr
		$\E^{- 4 x^2} (9+72 x^2 -192 x^4 + 512 x^6)$ 	 	&---		&  0.15 		&   2.87 		& 0.24  \cr
		$2 + {\cos(100 x)}$ 							&---		&  0.21		&  ---	 		& 0.11 \cr
		$\sech(200 x)$ 								&---		&  0.67		& 6.76	 	& 6.89  \cr
		\hline
 $ \E^{-100(x-1)^2} + \E^{-100(y+1)^2}(1+\cos(20x)) $		& 2  		& 1.34 		& --- 			& 2.46 \cr
		$\E^{-{x^4 \over 2} - {y^4 \over 2}} (x-y)^2$ 		&3 		& 0.54 		& ---  	  	& 6.41 	\cr
		$\E^{-x^2-2y^2}  \sech(10 x y)$ 					& 16 		& 10.90 		& ---  	  	& 9.43	\cr 
                $\E^{-x^2-2y^2}   (y-x)^2\sech(10 x y)$ 			& 51 		& 19.03 		& ---  	  	& 7.12 
\end{tabular}
\end{center}
\caption{\label{tab:SamplingTimes} Sampling  times  (in seconds) for our inverse transform sampling (ITS) compared to {\tt slicesample} in \MATLAB\ (SS), and rejection sampling (RS) using the function's maximum for  10,000 samples. The 2D inverse transform sampling approach described in Section \ref{sec:twodimensions} is remarkably efficient when the probability distribution is approximated by a low rank function.  (Note: {\tt slicesample} only works on decaying PDFs, hence fails on the third example.)}
\end{table}

\begin{figure}
	\includegraphics[width = .49\hsize]{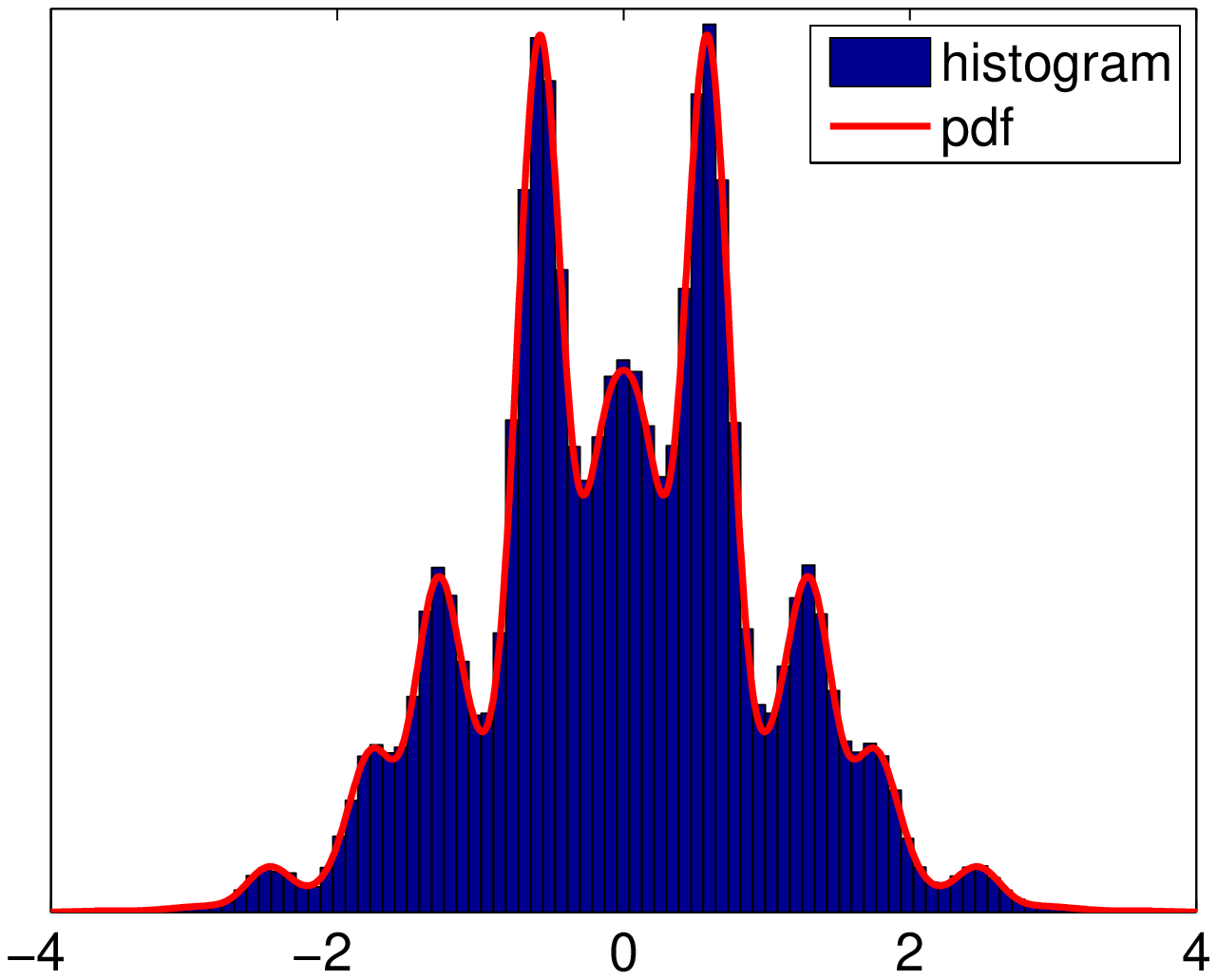}
	\includegraphics[width = .49\hsize]{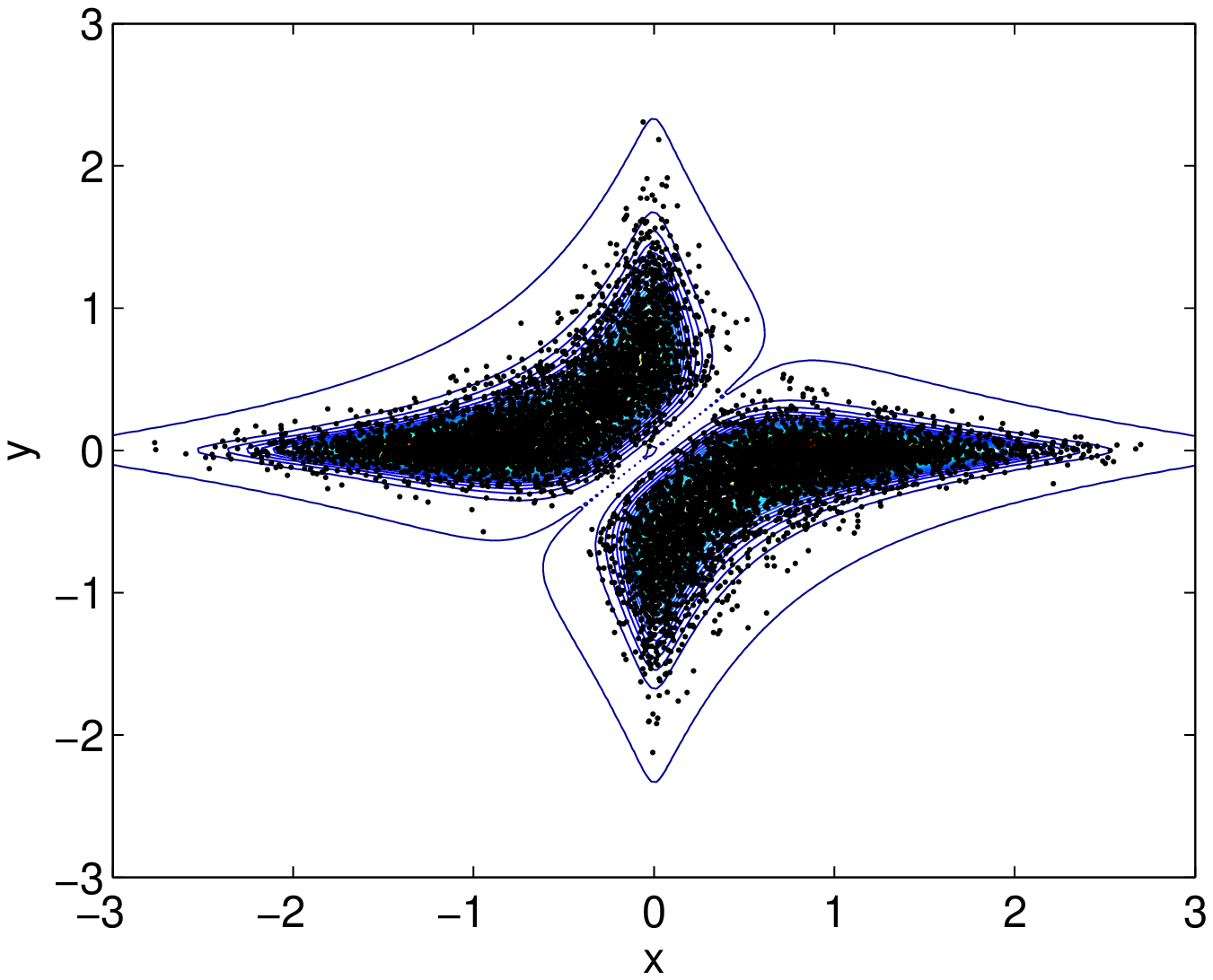}
\caption{\label{examples}
Left: histogram of one hundred thousand samples from the multimodal density computed in 2.28 seconds, and it can be seen that the generated pseudo-random samples have the correct distribution.  Right: ten thousand pseudo-random samples from the butterfly density: $f(x,y) = \E^{-x^2-2y^2} \sech(10 x y)(x-y)^2$. Only one 
sample lies out of the contour lines (the last contour is at the level curve $0.0001$), and this is consistent with the number of generated samples.   
 }
\end{figure}

In Table~\ref{tab:SamplingTimes}, we compare the sampling times\footnote{On a 2011 2.7 Ghz Intel Core i5 iMac.} of our algorithms based on inverse transform sampling, the slice sampling ({\tt slicesample}) routine in \MATLAB, and rejection sampling with a rectangular hat function bounded above by the exact maximum of the distribution.   
We consider the following univariate distributions:
\begin{enumerate}
	\item Multimodal density: $\E^{-{x^2 \over 2}} (1+ \sin^2 (3 x))(1+ \cos^2(5 x))$ on $[-8,8]$,
	\item The spectral density of a $4 \times 4$ Gaussian Unitary Ensemble \cite{DeiftOrthogonalPolynomials}: $\E^{- 4 x^2} (9+72 x^2 -192 x^4 + 512 x^6)$ on $[-4,4]$,
	\item Compactly supported oscillatory density: $2 + {\cos(100 x)}$	 on $[-1,1]$,
	\item Concentrated sech density: $\sech(200 x)$ on $[-1,1]$,
\end{enumerate}
and the following bivariate distributions:
\begin{enumerate}
\item Bimodal distribution: $\E^{-100(x-1)^2} + \E^{-100(y+1)^2}(1+\cos(20x))$ on $[-2,2]\times [-2,2]$,
	\item Quartic unitary ensemble eigenvalues \cite{DeiftOrthogonalPolynomials}: $\E^{-{x^4 \over 2} - {y^4 \over 2}} (x-y)^2$ on $[-7,7]\times [-7,7]$,
	\item Concentrated 2D sech density: $\E^{-x^2-2y^2}  \sech(10 x y)$ on $[-5,5] \times [-4,4]$,
	\item Butterfly density: $\E^{-x^2-2y^2}  \sech(10 x y) (x-y)^2$ on $[-3,3] \times [-3,3]$.
\end{enumerate}
Figure~\ref{examples} demonstrates our algorithm. Figure~\ref{examples} (left) shows a histogram of one hundred thousand pseudo-random samples generated from the multimodal density, and  Figure~\ref{examples} (right) shows ten thousand pseudo-random samples from the butterfly density. 

We observe that the inverse transform sampling approach described in Section~\ref{sec:inversesampling} significantly outperforms the \MATLAB\ implementation of slice sampling in 1D.   In the worst case, the computational cost is comparable to rejection sampling, and in two of 
the examples it outperforms rejection sampling by a factor of 10.

\begin{figure}
	\includegraphics[width = .49\hsize]{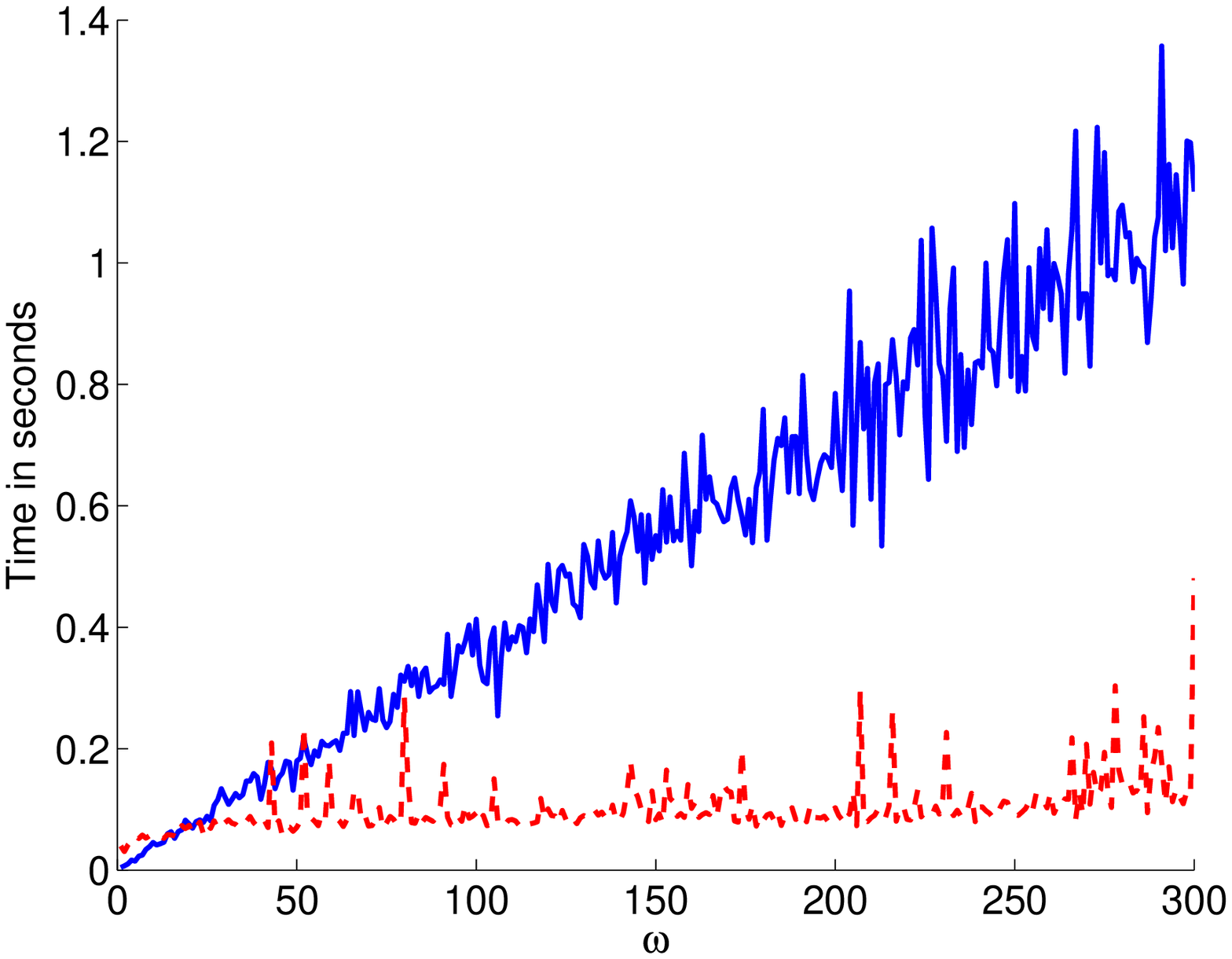}
	\includegraphics[width = .49\hsize]{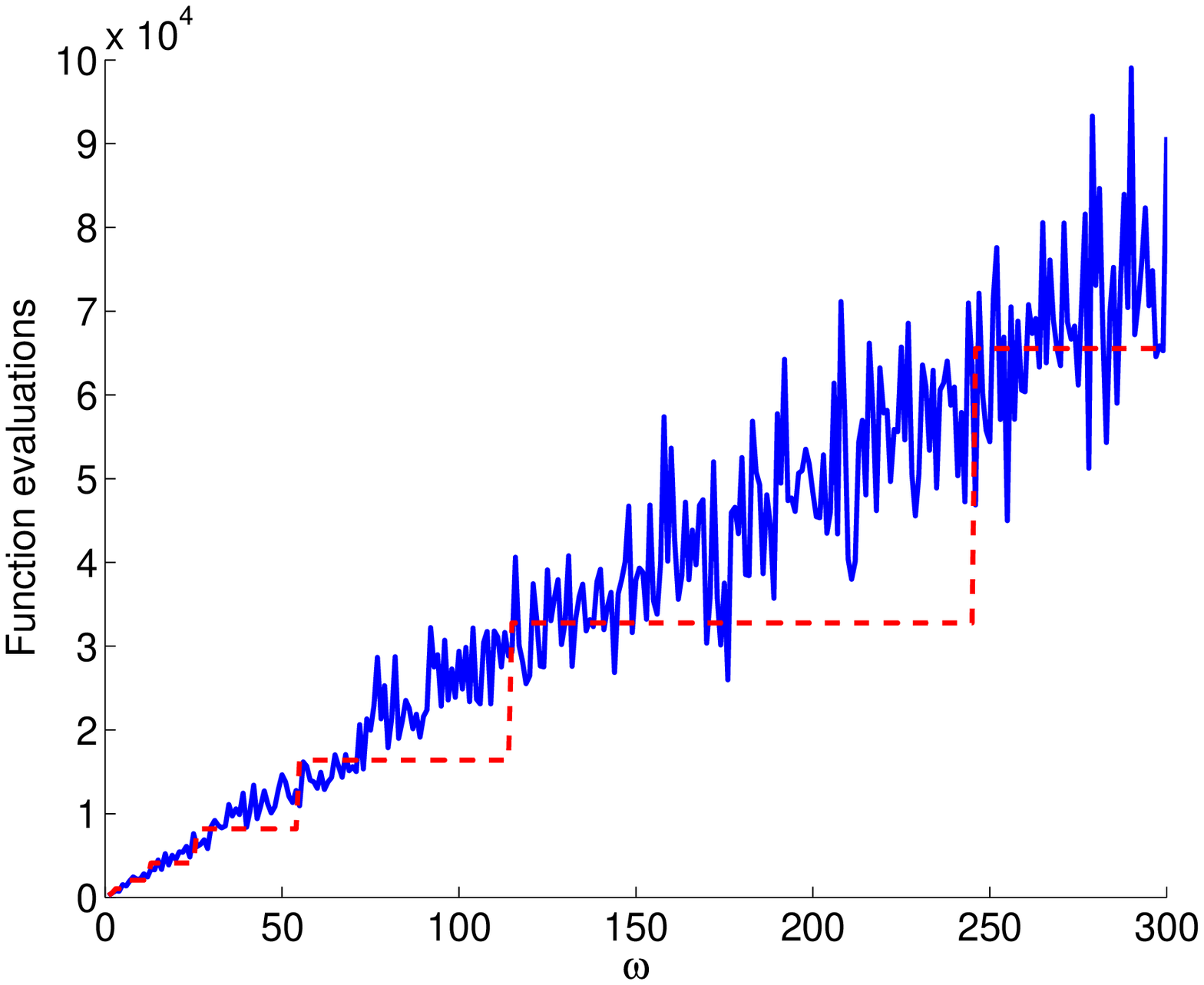}
\caption{\label{VsRejection}
Left: time (in seconds) of rejection sampling (solid) versus inverse transform method (dashed) for 100 samples of $\sech(k x)$.  Right: the number of function evaluations for 50 pseudo-random sample generated by rejection sampling (solid) versus the number of function evaluations for any number of inverse transform samples (dashed). 
 }
\end{figure}

Rejection sampling has practical limitations that inverse sampling avoids.  Finding the exact maximum of the distribution for the  rejection sampling is computationally expensive, and replacing the exact maximum with a less accurate upper bound causes more samples to be rejected, which reduces its efficiency.  Furthermore, the computational cost of each sample is not bounded, and to emphasize this we consider generating pseudo-random samples from $\sech(\omega x)$ for varying $\omega$ on the interval $[-8,8]$.   For large $\omega$ 
the percentage of the rectangular hat functions area that is under the probability distribution is small, and the rejection sampling approach rejects a significant proportion of its samples, while 
for our approach the dominating computational cost stays roughly constant.
In Figure~\ref{VsRejection} (left) we compare the computational cost of rejection sampling and the approach described in Section~\ref{sec:inversesampling}, and we observe that the Chebyshev expansion outperforms rejection sampling for $\omega\geq 30$.

\begin{remark}
	Of course, for a particular probability distribution a better hat function may be known, which would improve the efficiency of rejection sampling. However, in general accurate hat functions are not available and a rectangular hat function with a known upper bound is the practical choice.
\end{remark}


In many applications, the actual evaluation of the probability distribution is expensive. In such circumstances, another benefit of using polynomial approximation in the inverse transform sampling method is that 
once an approximation has been calculated the original distribution can be discarded, whereas rejection sampling evaluates the original distribution several 
times per sample.  Thus, in Figure~\ref{VsRejection}, we compare the number of function evaluations required to generate $50$ samples using rejection sampling versus the number of 
evaluations for our inverse transform sampling approach. The number of function evaluations are comparable at 50 samples, and hence if 500 pseudo-random samples were required then the rejection sampling approach would require about ten times the number of function evaluations.  

Parallelization of the algorithm is important for extremely large sample sizes.  This is easily accomplished using a polynomial approximation since the parameters used to represent the CDF can  be stored in local memory and the task of 
generating $N$ pseudo-random samples can be divided up among an arbitrary number of processors or machines. However, the rejection sampling approach cannot be 
parallelized with a black box distribution because evaluation would cause memory locking.


\section{Conclusions}\label{sec:conc}

We have shown that inverse transform sampling is efficient when employed with Chebyshev polynomial approximation in one dimension, and with low rank function approximation in two dimensions. 
This allows for an automated and robust algorithm for generating pseudo-random samples from a large class of probability distributions.  Furthermore, the approach that we have described can be generalized to higher dimensions by exploiting low rank approximations of tensors.
	
In fact, the inverse transform sampling approach is efficient for any class of probability distributions that can be numerically approximated, evaluated, and integrated. Therefore, there are straightforward extensions to several other classes of distributions of one variable:
\begin{enumerate}
	\item Piecewise smooth probability distributions, utilizing piecewise polynomial approximation schemes and automatic edge detection \cite{Pachon_10_01},
	\item Probability distribution functions with algebraic endpoint singularities, 
	\item Probability distributions with only algebraic decay at $\pm \infty$, by mapping to a bounded interval via ${1 - x \over 1 + x}$.
\end{enumerate}
The software package {\sc Chebfun} \cite{Chebfun} implements each of these special cases.  

In two dimensions and higher, approximating general piecewise smooth functions is a more difficult problem, and many fundamental algorithmic issues remain.  
However, if these constructive approximation questions are resolved the inverse transform sampling approach will be immediately applicable.


\end{document}